\documentclass[12pt]{article}
\usepackage{graphicx}
\usepackage{color}
\usepackage{amsmath,amsxtra,amssymb,latexsym, amscd}
\usepackage{amsmath, amssymb}
\usepackage{amsthm}
\usepackage[left=3cm,right=2.5cm,top=2.5cm,bottom=2.5cm]{geometry}

\newtheorem{thm}{Theorem}

\newtheorem{lem}[thm]{Lemma}

\theoremstyle{plain}
\theoremstyle{definition}
\newtheorem{rem}[thm]{Remark}

\DeclareMathOperator{\divv}{div}

\begin{document}
\newcommand{\lon}{\longrightarrow}
\newcommand{\lom}{\longmapsto}
\newcommand{\dis}{\displaystyle}

\title{\bf The classification of constant weighted curvature curves in the plane with a log-linear density}

\author{\bf  Doan The Hieu\\
 Departement of Mathematics\\
 College of Education, Hue University, Hue, Vietnam\\
 dthehieu@yahoo.com\\
{\bf Tran Le Nam} \\
 Departement of Mathematics\\
 Dong Thap University, Dong Thap, Vietnam\\
 lenamdongthapmuoi@gmail.com}

\maketitle
\begin{abstract}
 In this paper, we classify the class  of constant weighted curvature curves in the plane with a log-linear density, or in other words,  classify all traveling curved fronts with a constant forcing term in $\Bbb R^2.$ The classification gives some interesting phenomena and consequences including: the family of curves  converge to a round point when the weighted curvature of curves (or equivalently the forcing term of traveling curved fronts) goes to infinity, a simple proof for a main result in \cite{Jian2} as well as some well-known facts concerning to the isoperimetric problem in the plane with density $e^y.$

 \end{abstract}
\vskip 1cm
\section{Introduction}

Consider the following equation
\begin{equation}\label{equ1}\divv\left(\frac{\nabla u}{\sqrt{1+|\nabla u|^2}}\right)=\frac{a}{\sqrt{1+|\nabla u|^2}}+b,\end{equation}
where $u:U\subset\Bbb R^n\longrightarrow \Bbb R,$\ $U$ is an open domain  and $b$ is a constant.
This is the equation for self-similar translators with constant velocity $a$ of a mean curvature flow (MCF) with constant forcing term ${\cal H}=b$
\begin{equation}\label{equ2}\frac{\partial u}{\partial t}=\sqrt{1+|\nabla u|^2}\divv\left(\frac{\nabla u}{\sqrt {1+\nabla u|^2}}-{\cal H}\right).\end{equation}
MCF (\ref{equ2}) can be thought of as a generalization of the classical MCF when ${\cal H}=0$ and $a=1$ or as a modified version of MCF with constant forcing term in Minkowski space, a problem studied by Ecker and Huisken \cite{EcHui1}, followed by Aarons \cite{Aa}, H. Jian et al. \cite{Jian1}, \cite{Jian3} and others.

Equation (\ref{equ1}) appears in the study of traveling curved fronts, originated by the research of vortex motion of Ginzburg-Landau superconduction equation or by the research of interfacial phenomena,  and in several other fields in Physics (see \cite{NiTa}, \cite{Jian2}).

It also appears in the study of constant weighted mean curvature hypersurfaces in $\Bbb R^{n+1}$ with a log-linear density.

 The study of the geometry of manifolds with density has increased in the last seven years after the appearance of Morgan's paper ``Manifolds with density" published in the Notices Amer. Math. Soc. journal in 2005. A manifold with density is a Riemannian manifold with a positive function, called the density, used to weight both volume and hypersurface area. On a manifold with density $f=e^{\varphi},$ the weighted mean curvature of a hypersurface with unit normal $N$ is defined as follows
\begin{equation} \label{Hgro}   H_f:=H-\frac 1{n}\frac{d\varphi}{d{N}},\end{equation}
where $H$\ is the (Riemannian) mean curvature  of the hypersurface.

 For more details about manifolds with density and some relative topics we refer the reader to \cite{caja}-\cite{dadu}, \cite{Hiho}, \cite{Hi}, \cite{mamo}-\cite{mo5}, \cite{RCBM}.

 A typical example of manifolds with density is  Gauss space $G^n,$ the Euclidean space $\Bbb R^n$ with Gaussian probability density $f=(2\pi)^{-\frac n2}e^{-\frac{r^2}2},$ that is very interesting to probabilists. Weighted minimal hypersurfaces, i. e. $H_f=0,$  in Gauss spaces are known as self-similar shrinkers. This is one of three types of self-similar solutions of MCF  (self-similar shrinkers, self-similar translators and self-similar expanders) that are most studied recently.
  Self-similar translators are just weighted minimal hypersurfaces in $\Bbb R^{n+1}$ with a log-linear density $f,$ that is  $f=e^\varphi,$ where $\varphi(x_1, x_2, \ldots, x_{n+1})= \sum_{i=1}^{n+1} a_ix_i+a_0$ (see \cite{Hiho}).

Indeed, let $\Gamma_u\subset \Bbb R^{n+1}$ be the graph of the function $u.$ Because $\divv\left({\nabla u}/{\sqrt{1+|\nabla u|^2}}\right)=nH,$  and $a/{\sqrt{1+|\nabla u|^2}}=a\langle N, e_n\rangle,$ where $H$ is the mean curvature of $\Gamma_u$ and $N$ is the unit normal of the graph; Equation (\ref{equ1}) can be rewritten as
\begin{equation}\label{equ3} H_f=c,\end{equation}
where $c=b/n$ and $H_f:=H-(a/n)\langle N, e_{n+1}\rangle$ is the weighted mean curvature of $\Gamma_u$ with density $e^{ax_{n+1}}.$

Therefore, solutions of  Equation (\ref{equ1}) are constant weighted mean curvature hypersurfaces in $\Bbb R^{n+1}$ with a log-linear density $f=e^{ax_{n+1}}.$ Weighted minimal hypersurfaces in $\Bbb R^{n+1}$ with the density $f=e^{x_{n+1}}$ (self-similar translators of MCF) are solutions of (\ref{equ3}) for $a=1$ and $c=0.$

The classification of all constant weighted curvature hypersurfaces as well as weighted minimal hypersurfaces  in $\Bbb R^{n+1}$ with a general log-linear density is equivalent to the one  with the density $f=e^{x_{n+1}}$ (see Section \ref{secdensity}). In other words, we can assume that the velocity $a=1$ in the  classification of all self-similar translators of the MCF (\ref{equ2}).

Self-similar translators of MCF play a crucial role in the study type II singularities of the flow (see \cite{HuSi1}, \cite{HuSi2}). Beside lines, the first example of self-similar translators in the plane, the Grim Reaper, was discovered by Calabi. For higher dimensional cases, only a few specific examples of self-similar translators are known. In terms of weighted minimal surfaces in $\Bbb R^3$ with density $e^z,$ some cylindrical type examples have been showed in \cite{Hiho}.

In this paper, we solve Equation (\ref{equ3}) for parametric curves in the simplest case $n=1$ and therefore give the classification of constant weighted curvature curves in the plane with a log-linear density, or in other words,  the classification of all traveling curved fronts with a constant forcing term in $\Bbb R^2.$ Although the classification is elementary, it gives some interesting phenomena and consequences including: the family of curves  converge to a round point when the weighted curvature of curves goes to infinity (see Section \ref{convergence}); a simple proof for a main result in \cite{Jian2} as well as some well-known facts concerning to the isoperimetric problem in the plane with density $e^y.$

{\bf Acknowledgement.}\ We would like to thank Professor Frank Morgan for his encouragement and for reading the first draft of the paper. The authors are supported in part by the National Foundation for Science and Technology Development, Vietnam (Grant No. 101.01-2011.26).

\section{Spaces with log-linear densities}\label{secdensity}
Consider $\Bbb R^{n}$ with a log-linear density $f=e^\varphi,$ where $\varphi(x_1, \ldots, x_{n+1})= \sum_{i=1}^{n+1} a_ix_i+a_0.$ In this space, the weighted mean curvature of a hypersurface with the unit normal  $N$ is
$$H_f:=H-\frac1{n}\langle N, e_{n}\rangle,$$

The quantity $\langle N, e_{n+1}\rangle$ is just the projection of $N$ on the $x_{n+1}$-axis. This quantity and the mean curvature are invariant under rigid motions. Therefore,  when study any problem concerning only to weighted mean curvature, without loss of generality we can assume that $\varphi = ax_{n+1}.$ Moreover, we can assume $a=1.$

Now consider the plane with density $e^{ay}.$ This is a  plane of zero generalized Gauss curvature, $G_\varphi=G-\triangle\varphi,$ where $G$ is the classical Gauss curvature. The notion of generalized Gauss curvature was introduced in Corwin {\sl et al.} \cite{co2} to study the generalization of the Gauss-Bonnet formula for surfaces with densities, which is finally found in the most general form for the case of unrelated densities  in \cite{co3}.

Let $\alpha(s)=(x(s),y(s))$ be a parametric curve with arc length parameter $s,$ and let $\beta(s)=(1/a)(x(as), y(as)).$ Denote by $k_f(\alpha)$ the weighted curvature of  $\alpha$ corresponding to the density $f=e^{y}$ and $k_{\widetilde{f}}(\beta)$ the weighted curvature of  $\beta$ corresponding to the density $\widetilde{f}=e^{ay}.$ By a simple computation, it is showed that the weighted curvature of  $\alpha$ in $\Bbb R^2$ with density $f=e^{y}$ is
\begin{equation}\label{wc}
  k_f(\alpha)=x'y''-x''y'-x',
\end{equation}
 while the weighted curvature of  $\beta$ in $\Bbb R^2$ with density $\widetilde{f}=e^{ay}$ is
\begin{equation}\label{wc2}
  k_{\widetilde{f}}(\beta)=ax'y''-ax''y'-ax'.
\end{equation}
Therefore, we get
\begin{lem} \label{wc11} $$ak_f(\alpha)=k_{\widetilde{f}}(\beta).$$
\end{lem}

\section{Constant weighted curvature curves in the plane with density $e^y$}
For solving Equation (\ref{equ3}) to classify all constant weighted curvature curves in the plane with density $e^{ay},$ by Lemma \ref{wc11}, we can assume that $a=1.$ Moreover, because the weighted curvature of a curve is invariant under translations,  some constant of integration  will be omitted for simplicity.

Equation (\ref{equ3}) for a parametric curve $\alpha(s)=(x(s),y(s))$  with arc length parameter is
\begin{equation}\label{wc111}
  x'y''-x''y'-x'=c.
\end{equation}

Set
\begin{align}
\label{pt:thamso}
\begin{cases}
    x'=-\cos(2\xi)\\
    y'=\sin(2\xi)
\end{cases}.
\end{align}

Then, the equation (\ref{wc111}) reads as follows:
\begin{align}\label{pt:dchang}
 -2\xi'+\cos(2\xi)=c.
\end{align}

If $\cos (2\xi)\ne c$ for all  $s,$ then the equation (\ref{pt:dchang}) is rewritten as follows
\begin{align}
  \label{pt:docong1}
-2\xi'=\dfrac{c-1+(c+1)\tan^2\xi}{1+\tan^2\xi}.
\end{align}


\subsection{The case of $\boldsymbol{c<-1}$ (see Figure 6.1)}

Equation  (\ref{pt:docong1}) can be rewritten as
\begin{align*}
    &-\dfrac{2d\tan\xi}{\tan^2\xi+a^2}=(c+1)ds, \;\mbox{where } \, a=\sqrt{\dfrac{c-1}{c+1}}.
\end{align*}
Thus,
   $$ -\dfrac{2}{a}\arctan\dfrac{\tan\xi}{a}=(c+1)s$$
or
\begin{equation}\label{ccc} \tan\xi=a\tan\Big(\dfrac{\sqrt{c^2-1}}{2}s\Big),\end{equation}
because $c+1<0.$
From (\ref{ccc}), we have

  $$  \begin{cases}
      x'(s)=
      \dfrac{a^2\tan^2\Big(\dfrac{\sqrt{c^2-1}}{2}s\Big)-1}{1+a^2\tan^2\Big(\dfrac{\sqrt{c^2-1}}{2}s\Big)},\\
       y'(s)=
      \dfrac{2a\tan\Big(\dfrac{\sqrt{c^2-1}}{2}s\Big)}{1+a^2\tan^2\Big(\dfrac{\sqrt{c^2-1}}{2}s\Big)},\\
    \end{cases}$$
A direct computation yields
\begin{align}    \begin{cases}\label{tam}
          x'(s)=\dfrac{-2}{\sqrt{c^2-1}} \dfrac{\Big(1-a^2\tan^2\dfrac{\sqrt{c^2-1}}{2}s\Big)d\Big(\tan\dfrac{\sqrt{c^2-1}}{2}s\Big)}
      {\Big(1+a^2\tan^2\dfrac{\sqrt{c^2-1}}{2}s\Big)\Big(1+\tan^2\dfrac{\sqrt{c^2-1}}{2}s\Big)},\\[0.8cm]
      y'(s)=\dfrac{-4}{(c+1)}\dfrac{\tan\Big(\dfrac{\sqrt{c^2-1}}{2}s\Big)d\Big(\tan\dfrac{\sqrt{c^2-1}}{2}s\Big)}
      {\Big(1+a^2\tan^2\dfrac{\sqrt{c^2-1}}{2}s\Big)\Big(1+\tan^2\dfrac{\sqrt{c^2-1}}{2}s\Big)}.
    \end{cases} \end{align}

Integrating both sides of (\ref{tam}), we obtain

\begin{align}
    \begin{cases}\label{eq:gt1}
        x(s)=\dfrac{2(a^2+1)}{\sqrt{c^2-1}(a^2-1)}\dfrac{\sqrt{c^2-1}}{2}s-\dfrac{4a}{\sqrt{c^2-1}(a^2-1)}
        \arctan\left(\sqrt{\dfrac{c-1}{c+1}}\tan\dfrac{\sqrt{c^2-1}}{2}s\right),        \\[0.4cm]
        y(s)=\dfrac{2}{(c+1)(a^2-1)}\ln\left(\dfrac{\tan^2\left(\dfrac{\sqrt{c^2-1}}{2}s\right)+1}
        {\dfrac{c-1}{c+1}\tan^2\left(\dfrac{\sqrt{c^2-1}}{2}s\right)+1}\right).
    \end{cases}
\end{align}
\indent Substitute $a=\sqrt{\dfrac{c-1}{c+1}}$ into (\ref{eq:gt1}), we get
\begin{align*}
    \begin{cases}
        x(s)=-2\arctan\left(\sqrt{\dfrac{c-1}{c+1}}\tan\dfrac{\sqrt{c^2-1}}{2}s\right)-cs,\\[0.4cm]
        y(s)=-\ln\left(\dfrac{\tan^2\left(\dfrac{\sqrt{c^2-1}}{2}s\right)+1}
        {\dfrac{c-1}{c+1}\tan^2\left(\dfrac{\sqrt{c^2-1}}{2}s\right)+1}\right).
    \end{cases}
\end{align*}
\subsection{The case of $\boldsymbol{c=-1}$ (see Figure 6.2)}
If there exists  $s_0$ so that $\cos[2\xi(s_0)]=-1,$ then $\xi(s)=\xi(s_0)$ is the unique solution of the equation (\ref{pt:dchang}) and the corresponding curves are straight lines parallel to the $x$-axis.

If $\cos 2\xi(s)\neq -1$ for all  $s,$ then the equation (\ref{pt:docong1})  becomes
$$\xi'=\dfrac{1}{1+\tan^2\xi}.$$
Thus,
$$d(\tan\xi)=ds,$$
and therefore,
$$   \tan\xi = s.$$
 Then, by Equation (\ref{pt:thamso})
\begin{align}\label{123}
    \begin{cases}
    x'(s)=-\dfrac{1-s^2}{1+s^2},\\[0.2cm]
    y'(s)=\dfrac{2s}{1+s^2}.
  \end{cases}
\end{align}
Integrating both sides (\ref{123}), we obtain
\begin{align} \label{cbang-1}\begin{cases}
      x(s)=-2\arctan s+s,\\
    y(s)=\ln(1+s^2).
  \end{cases}\end{align}
\subsection{The case of $\boldsymbol{-1 < c <1}$ (see Figures 6.3, 6.4, 6.5)}
If there  exists  $s_0$ so that $\cos[2\xi(s_0)]=c,$ then $\xi(s)=\xi(s_0)$ is the unique solution of the equation (\ref{pt:dchang}) and the corresponding  curves are straight lines with the slope $-\frac{\sqrt{1-c^2}}{c}.$ In the case of $c=0,$ the slope of the lines is $\infty,$ i.e. the lines are vertical.

If this is not to be the case, Equation (\ref{pt:docong1}) is equivalent to
\begin{align}\label{eee}
    &-\dfrac{2d(\tan\xi)}{\tan^2\xi-b^2}=(c+1)ds, \;\mbox{where } \, b=\sqrt{\dfrac{1-c}{c+1}}\neq0.\end{align}
Solving (\ref{eee}) yields
  $$  \dfrac{1}{b}\ln\Big|\dfrac{\tan\xi+b}{\tan\xi-b}\Big|=(c+1)s,$$
or
$$    \Big|\dfrac{\tan\xi+b}{\tan\xi-b}\Big|={e}^{\sqrt{1-c^2}s}.$$
We have two cases to consider:

{\bf Case 1:}\ $ \dfrac{\tan\xi+b}{\tan\xi-b}={e}^{\sqrt{1-c^2}s},$ or equivalently,
$   \tan\xi=b\dfrac{{e}^{\sqrt{1-c^2}s}+1}{{e}^{\sqrt{1-c^2}s}-1}.$

In this case, we have

$$\begin{cases}
  x'(s)=-\dfrac{1-\Big(b\dfrac{{e}^{\sqrt{1-c^2}s}+1}{{e}^{\sqrt{1-c^2}s}-1}\Big)^2}
        {1+\Big(b\dfrac{{e}^{\sqrt{1-c^2}s}+1}{{e}^{\sqrt{1-c^2}s}-1}\Big)^2},\\
  y'(s)=\dfrac{2\Big(b\dfrac{{e}^{\sqrt{1-c^2}s}+1}{{e}^{\sqrt{1-c^2}s}-1}\Big)}
  {1+\Big(b\dfrac{{e}^{\sqrt{1-c^2}s}+1}{{e}^{\sqrt{1-c^2}s}-1}\Big)^2};
\end{cases}$$
or
\begin{equation}\label{111} \begin{cases}
 x'(s)=-\dfrac{(1-b^2){e}^{2\sqrt{1-c^2}s}-2(b^2+1){e}^{\sqrt{1-c^2}s}+1-b^2}
  {(1+b^2){e}^{2\sqrt{1-c^2}s}+2(b^2-1){e}^{\sqrt{1-c^2}s}+b^2+1},\\[0.2cm]
 y'(s)=\dfrac{2b\big({e}^{2\sqrt{1-c^2}s} - 1\big)}
  {(1+b^2){e}^{2\sqrt{1-c^2}s}+2(b^2-1){e}^{\sqrt{1-c^2}s}+b^2+1}.
\end{cases}\end{equation}

Substitute $b=\sqrt{\dfrac{1-c}{1+c}}$ into (\ref{111}), we get
$$\begin{cases}
      x'(s)=-\dfrac{\frac{2c}{c+1}{e}^{2\sqrt{1-c^2}s}-2\frac{2}{c+1}{e}^{\sqrt{1-c^2}s}+\frac{2c}{c+1}}
      {\frac{2}{c+1}{e}^{2\sqrt{1-c^2}s}+2\frac{-2c}{c+1}{e}^{\sqrt{1-c^2}s}+\frac{2}{c+1}}, \\[0.2cm]
        y'(s)=\dfrac{2\sqrt{\frac{1-c}{1+c}}\big({e}^{2\sqrt{1-c^2}s}-1\big)}
      {\frac{2}{c+1}{e}^{2\sqrt{1-c^2}s}+2\frac{-2c}{c+1}{e}^{\sqrt{1-c^2}s}+\frac{2}{c+1}};
\end{cases}$$
or
\begin{align} \label{fff}
\begin{cases}
  x'(s) = -\dfrac{c{e}^{2\sqrt{1-c^2}s} - 2{e}^{\sqrt{1-c^2}s} + c}
  {{e}^{2\sqrt{1-c^2}s} - 2c {e}^{\sqrt{1-c^2}s} + 1},\\[0.4cm]
   y'(s) = \dfrac{\sqrt{1-c^2}({e}^{\sqrt{1-c^2}s}-{e}^{-\sqrt{1-c^2}s})}
  {{e}^{\sqrt{1-c^2}s} - 2c + {e}^{-\sqrt{1-c^2}s}}.
\end{cases}
\end{align}

Integrating both sides of  (\ref{fff}), we obtain
\begin{equation}
\label{eq:nho1}
     \begin{cases}
          x(t)=2\,\arctan \left( {\dfrac {{{e}^{\sqrt {1-{c}^{2}}s}}-c}
          {\sqrt {1-{c}^{2}}}} \right)-cs, \\[0.4cm]
          y(s)=\ln  \left( {{e}^{\sqrt{1-c^2}s}  + {e}^{-\sqrt{1-c^2}s}} - 2c\right).
    \end{cases}
\end{equation}

{\bf Case 2:}\ $ \dfrac{\tan\xi+b}{\tan\xi-b}={e}^{\sqrt{1-c^2}s},$ or equivalently,
 $ \tan\xi=b\dfrac{{e}^{\sqrt{1-c^2}s}-1}{{e}^{\sqrt{1-c^2}s}+1}.$
After a quite similar computation as in Case 1, we obtain
\begin{equation}
\label{eq:nho1h}
\begin{cases}
            x(s)=-2\,\arctan \left( {\dfrac {{{e}^{\sqrt {1-{c}^{2}}s}}+c}
          {\sqrt {1-{c}^{2}}}} \right)-cs,\\[0.3cm]
          y(s)=\ln  \left( {{e}^{\sqrt{1-c^2}s}  + {e}^{-\sqrt{1-c^2}s}} + 2c\right).
    \end{cases}
\end{equation}

\begin{rem}\
\begin{enumerate}
\item
Set $f_1(c,s)=2\arctan \Big( {\dfrac {{{e}^{\sqrt {1-{c}^{2}}s}}-c}
          {\sqrt {1-{c}^{2}}}} \Big)-cs$, $f_2(c,s)=-2\,\arctan \Big( {\dfrac {{{e}^{\sqrt {1-{c}^{2}}s}}+c}
          {\sqrt {1-{c}^{2}}}} \Big)-cs$, $g_1(c,s)=\ln  \left( {{e}^{\sqrt{1-c^2}s}  + {e}^{-\sqrt{1-c^2}s}} - 2c\right)$ and
          $g_2(c,s)=\ln  \left( {{e}^{\sqrt{1-c^2}s}  + {e}^{-\sqrt{1-c^2}s}} + 2c\right),$ we can verify that
\begin{align*}
\begin{cases}
  f_1(c,s)=-f_2(-c,s),\\
  g_1(c,s)=g(-c,s).
\end{cases}
\end{align*}
Therefore the traces of the curves defined by equation (\ref{eq:nho1}) and (\ref{eq:nho1h}) are symmetric to the  $x$-axis.

\item  When $c=0,$ Equation (\ref{eq:nho1}) and (\ref{eq:nho1h}) are the equations of the Grim Reapers (see Figure 6.4)
\begin{equation}
  \begin{cases}   \label{khong2}
    x(s)=\pm 2\arctan(e^s),\\
    y(s)=\ln(e^s+e^{-s});
    \end{cases}\ \ \ \ \ s\in\Bbb R.
\end{equation}

\item From Equation (\ref{khong2}) we can deduce some simple facts (see also Figure 6.4).
\begin{enumerate}
\item  Two zero weighted curvature curves intersect in at most one point. Therefore, the geodesic connecting two points if exists is unique.
\item  If the difference of the $x$-coordinates of two given points is not less than $\pi,$ then there exists no geodesic connecting these points.
\item If the difference of the $x$-coordinates of two given points is less than $\pi,$ then there exists a unique zero weighted curvature curve connecting these points and this curve is  the shortest path which  can be proved directly.
\end{enumerate}
\end{enumerate}
\end{rem}

\color{black}
\subsection{The case of $\boldsymbol{c=1}$ (see Figure 6.2)}
If there exists  $s_0$ so that $\cos[2\xi(s_0)]=1,$ then $\xi(s)=\xi(s_0)$ is the unique solution of the equation (\ref{pt:dchang}) and the corresponding curves are straight lines parallel to the $x$-axis.

If $\cos 2\xi(s)\neq 1$ for all  $s,$ then the equation (\ref{pt:docong1})  becomes
$$-\xi'=\dfrac{\tan^2\xi}{1+\tan^2\xi} =\dfrac{1}{1+\cot^2s}.$$
Thus,
$$d(\cot\xi)=ds,$$
and therefore,
$$   \cot\xi = s.$$
 Then, by Equation (\ref{pt:thamso})
\begin{align}\label{1234}
    \begin{cases}
    x'(s)=\dfrac{1-s^2}{1+s^2},\\[0.2cm]
    y'(s)=\dfrac{2s}{1+s^2}.
  \end{cases}
\end{align}
Integrating both sides (\ref{1234}) yields
\begin{align} \label{cbang1}\begin{cases}
      x(s)=2\arctan s-s,\\
    y(s)=\ln(1+s^2).
  \end{cases}\end{align}

We can  see that the curves determined by (\ref{cbang-1}) and (\ref{cbang1}) have the same traces but have opposite directions.
\subsection{The case of $\boldsymbol{c>1}$ (see Figure 6.6)}

With similar arguments as in the case of $c<-1$, we have
\begin{equation}
  \tan\xi=-a\tan\Big(\dfrac{\sqrt{c^2-1}}{2}s\Big).
\end{equation}
Therefore,
\begin{align}
  \begin{cases}\label{tam2}
          x'(s)=\dfrac{-2}{\sqrt{c^2-1}} \dfrac{\Big(1-a^2\tan^2\dfrac{\sqrt{c^2-1}}{2}s\Big)d\Big(\tan\dfrac{\sqrt{c^2-1}}{2}s\Big)}
      {\Big(1+a^2\tan^2\dfrac{\sqrt{c^2-1}}{2}s\Big)\Big(1+\tan^2\dfrac{\sqrt{c^2-1}}{2}s\Big)},\\[0.8cm]
      y'(s)=\dfrac{-4}{(c+1)}\dfrac{\tan\Big(\dfrac{\sqrt{c^2-1}}{2}s\Big)d\Big(\tan\dfrac{\sqrt{c^2-1}}{2}s\Big)}
      {\Big(1+a^2\tan^2\dfrac{\sqrt{c^2-1}}{2}s\Big)\Big(1+\tan^2\dfrac{\sqrt{c^2-1}}{2}s\Big)}.
    \end{cases}
\end{align}
Integrating both sides of (\ref{tam2}) and substituting $a=\sqrt{\dfrac{c-1}{c+1}}$ into
the result, we get
\begin{align*}
    \begin{cases}
        x(s)=2\arctan\left(\sqrt{\dfrac{c-1}{c+1}}\tan\dfrac{\sqrt{c^2-1}}{2}s\right)-cs,\\[0.4cm]
        y(s)=-\ln\left(\dfrac{\tan^2\left(\dfrac{\sqrt{c^2-1}}{2}s\right)+1}
        {\dfrac{c-1}{c+1}\tan^2\left(\dfrac{\sqrt{c^2-1}}{2}s\right)+1}\right).
    \end{cases}
\end{align*}

\section{Classification of constant weighted curvature curves}
Combining the above results, up to translations, the classification of  constant weighted curvature curves in the plane with density $e^y$ is stated as follows.
\begin{thm}
 \begin{enumerate}
 \item A curve with weighted curvature zero is either a straight line (parallel to the $y$-axis) or the Grim Reaper defined by (see Figure 6.4)
  \begin{equation}
  \begin{cases}
     x(s)=2\arctan(e^s),\\
     y(s)=\ln(e^s+e^{-s}).
     \end{cases} s\in\mathbb{R}.
\end{equation}
\item A curve with constant weighted curvature $|k_\varphi|<1$ is either a straight line  or the one defined by (see Figure 6.3, Figure 6.5)

 $$\begin{cases}
          x(s)=2\,\arctan \left( {\dfrac {{{e}^{\sqrt {1-{c}^{2}}s}}- c}
          {\sqrt {1-{c}^{2}}}} \right) - cs,\\[0.4cm]
          y(s)=\ln\left( {{e}^{\sqrt{1-c^2}s}  + {e}^{-\sqrt{1-c^2}s}} - 2c\right).
    \end{cases} s\in\mathbb{R}.
$$

\item A curve with constant weighted curvature $\pm 1$ is either a straight line (parallel to the $x$-axis) or the one defined by (see Figure 6.2)

$$\begin{cases}
    x(s)=2\arctan s - s,\\
    y(s)=\ln(1+s^2).
  \end{cases} s\in\mathbb{R}.$$

  \item A curve with constant weighted curvature $|k_\varphi|>1$ is defined by  (see Figure 6.1, Figure 6.6)
  \begin{align*}
    \begin{cases}
             x(s)=\pm 2\arctan\left(\sqrt{\dfrac{c-1}{c+1}}\tan\dfrac{\sqrt{c^2-1}}{2}s\right)-cs,\\[0.4cm]
             y(s)=-\ln\left(\dfrac{\tan^2\left(\dfrac{\sqrt{c^2-1}}{2}s\right)+1}
        {\dfrac{c-1}{c+1}\tan^2\left(\dfrac{\sqrt{c^2-1}}{2}s\right)+1}\right).
    \end{cases}s\in\Big(-\dfrac{\pi}{\sqrt{c^2-1}},\dfrac{\pi}{\sqrt{c^2-1}}\Big).
\end{align*}
   \end{enumerate}
 \end{thm}
\section{Some consequenses}
\label{convergence}
\begin{enumerate}
\item The figures in the next section seem to show that if the weighted curvature $c$ goes to $\pm \infty$ then the limit of the curves is a point. For looking the behaviour of curves  near the limit point we can use the standard technique of rescaling the curves by a scale factor $\sqrt{c^2-1}.$ We give a proof for the case of $c\rightarrow\infty.$ The proof for $c\rightarrow-\infty$ is quite similar. Suppose that $\alpha (c)$ is a curve with arc length parameter and of constant weighted curvature $c>1.$ Let $\beta(c)=\sqrt{c^2-1}\alpha(c).$ The curvature of $\beta$ is $(1/\sqrt{c^2-1})k(\alpha(c)),$ where $k(\alpha(c))$ is the curvature of  the curve $\alpha (c).$ By Equation (\ref{wc111}), $k(\alpha (c))=x'+c.$ Therefore
    $$ \frac 1{\sqrt{c^2-1}}k(\alpha(c))=\frac 1{\sqrt{c^2-1}}\left(\dfrac{-c\cos(\sqrt{c^2-1}s)-1}{c+\cos(\sqrt{c^2-1}s)} +c\right).$$
 It is not hard to check that
 $$\lim_{c\rightarrow\infty}\Big|\frac 1{\sqrt{c^2-1}}k(\alpha(c))\Big| =1.$$
Thus, the family of curves converges to a round point when $c$ goes to infinity.
 \item The study of traveling fronts of curve flow with external force field for the simple case $\nabla w=(c_1,c_2)$ leads to the following equation (see \cite{Jian2}):
\begin{equation}\label{equJian} c=\frac{\varphi''(x)}{1+\varphi(x)^2}+c_2-c_1\varphi'(x).\end{equation}
One of the main results in \cite{Jian2} stated that, the solution of (\ref{equJian}) with an initial conditions is either  a line or the Grim Reaper.
In terms of weighted curvature, solutions of  (\ref{equJian}) are just zero weighted curvature curves in the plane with density $f=e^{-c_1x+(c_2-c)y},$ which are also zero weighted curvature curves in the plane with density $f=e^y$ under a suitable change of coordinates. In other words, the solutions are self-similar translators in the plane which are known to be either lines or Grim Reapers.
\item Proposition 4.8 in \cite{caja}, which states {\sl ``The plane with density $e^x$ contains no isoperimetric region''}, can be deduced from the classification. It is clear from the classification (and the figures) that  isoperimetric curves, i.e. curves bounding isoperimetric regions (in case of existence), must have infinite weighted length or singularities.
\end{enumerate}

  \section{Figures of constant weighted curvature curves}    \label{sec4}

\begin{tabular}{p{0.55\textwidth}c}
  \begin{minipage}[!htb]{0.5\textwidth}
     \centering
    \includegraphics[width=0.85\textwidth]{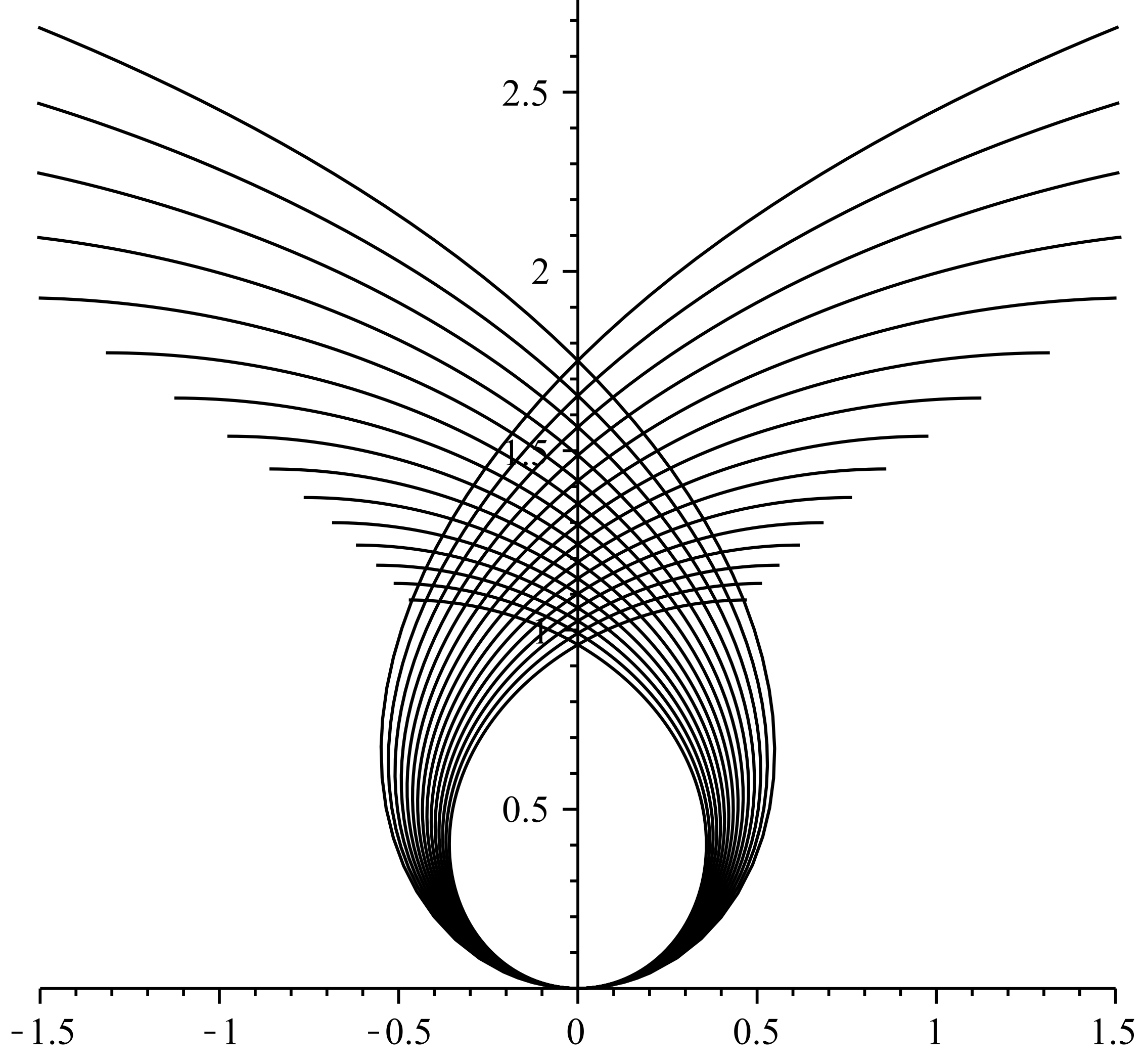}\\
    {Figure 6.1. Curves of $k_\varphi<-1$}.
 \end{minipage}
&
\begin{minipage}[!htb]{0.5\textwidth}\centering
    \includegraphics[width=0.85\textwidth]{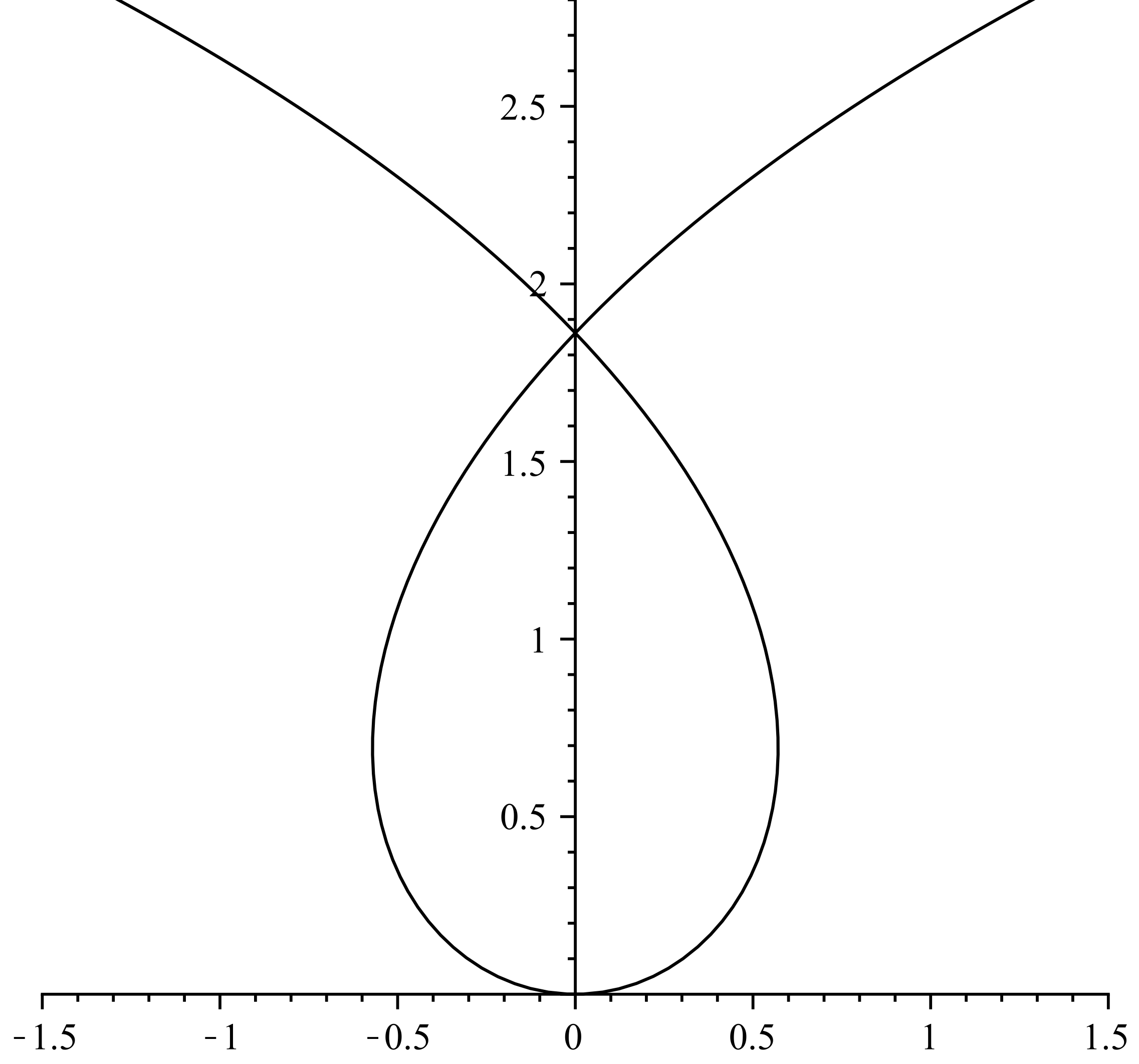}\\
    {Figure 6.2. Curves of  $k_\varphi=\pm1$}.
\end{minipage}\\
$\begin{cases}
             x(s)=- 2\arctan\left(\sqrt{\dfrac{c-1}{c+1}}\tan\dfrac{\sqrt{c^2-1}}{2}s\right)-cs,\\[0.4cm]
             y(s)=-\ln\left(\dfrac{\tan^2\left(\dfrac{\sqrt{c^2-1}}{2}s\right)+1}
        {\dfrac{c-1}{c+1}\tan^2\left(\dfrac{\sqrt{c^2-1}}{2}s\right)+1}\right).
     \end{cases}$&
$\begin{cases}
     x(s)=s-2\arctan s,\\
     y(s)=\ln(1+s^2).
  \end{cases}$
 \end{tabular}
\vskip0.5cm
\noindent
\begin{tabular}{cc}
\begin{minipage}[!thb]{0.5\textwidth}
\centering
 \includegraphics[width=0.85\textwidth]{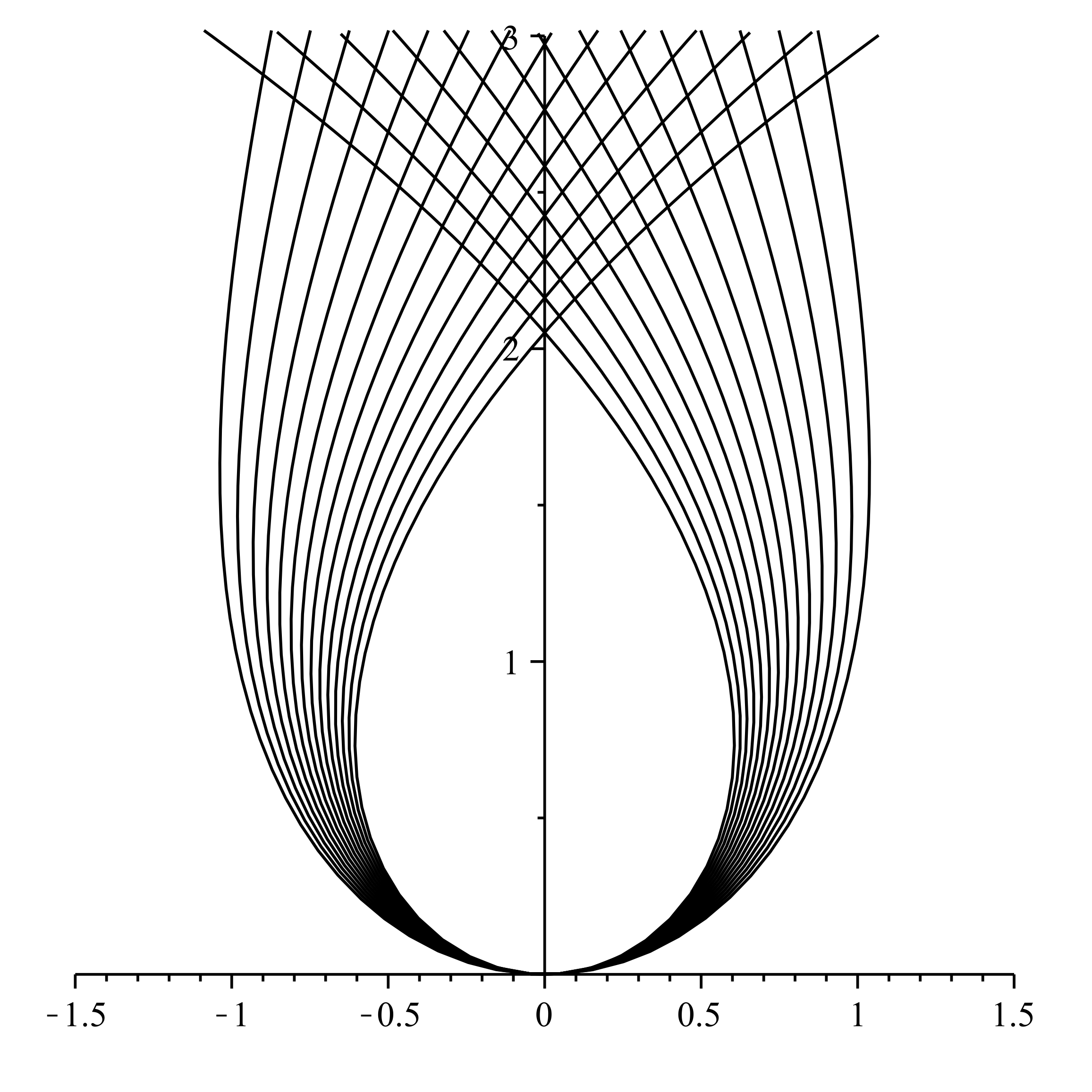}\\
 {Figure 6.3. Curves of  $k_\varphi\in(-1,0)$.}
 \end{minipage}
&\begin{minipage}[!thb]{0.5\textwidth}\centering
 \includegraphics[width=0.85\textwidth]{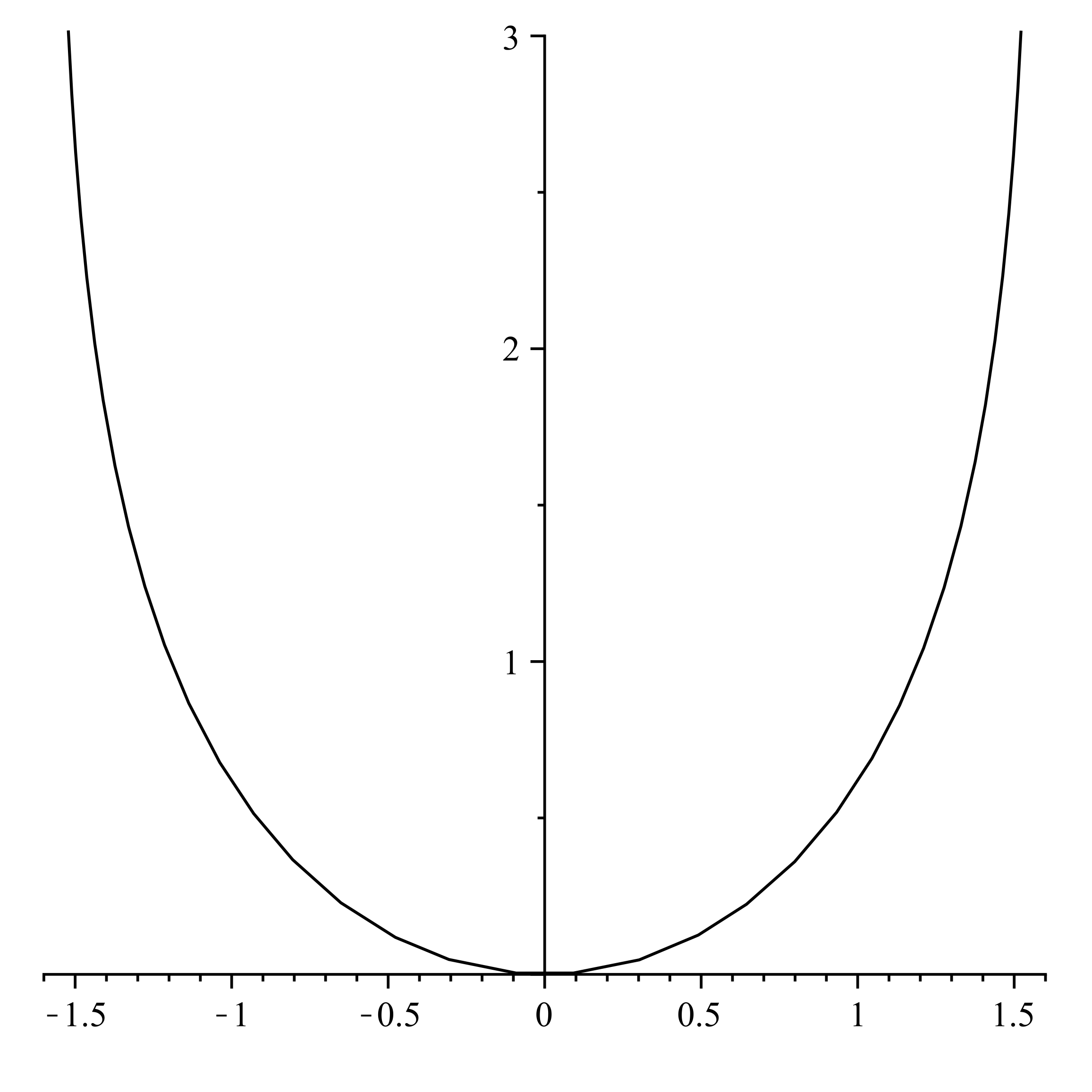}\\
 {Figure 6.4. Curve of  $k_\varphi=0$}.
\end{minipage}\\
   $\begin{cases}
          x(s)=2\,\arctan \left( {\dfrac {{{e}^{\sqrt {1-{c}^{2}}s}}-c}
          {\sqrt {1-{c}^{2}}}} \right) - cs,\\[0.4cm]
          y(s)=\ln\left( {{e}^{\sqrt{1-c^2}s}  + {e}^{-\sqrt{1-c^2}s}}- 2c\right).
    \end{cases}$&
    \begin{minipage}[!thb]{0.5\textwidth}
     \centering
     (the Grim Reaper curve)\\
    $\begin{cases}
     x(s)=2\arctan(e^s),\\
     y(s)=\ln(e^s+e^{-s}).
   \end{cases}$
   \end{minipage}
\end{tabular}
\begin{tabular}{cc}
\begin{minipage}[!thb]{0.5\textwidth}
 \centering
 \includegraphics[width=0.9\textwidth]{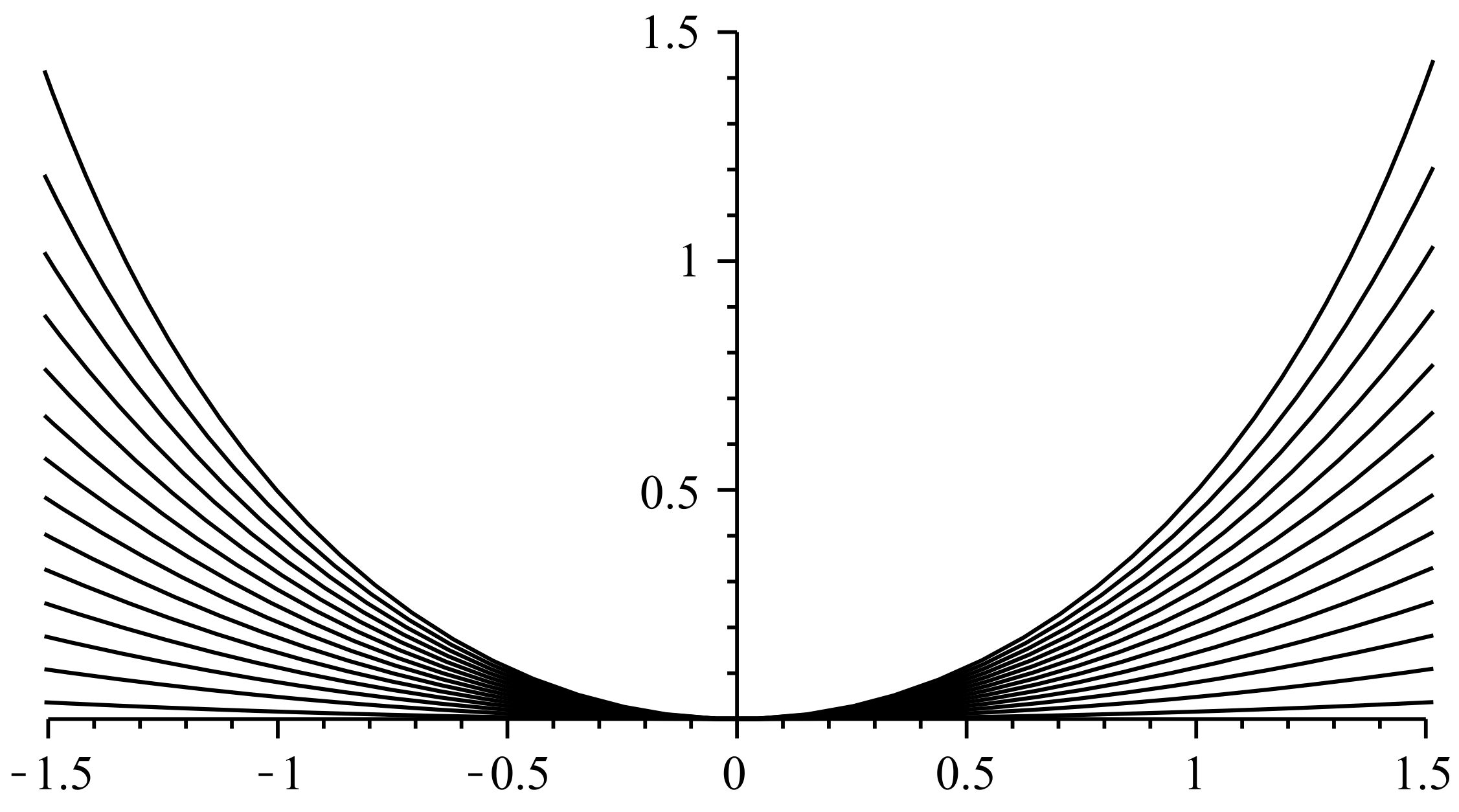}\\
 {Figure 6.5. Curves of $k_\varphi\in(0,1)$}.
\end{minipage}&
\begin{minipage}[!thb]{0.5\textwidth}
 \centering
  \includegraphics[width=0.9\textwidth]{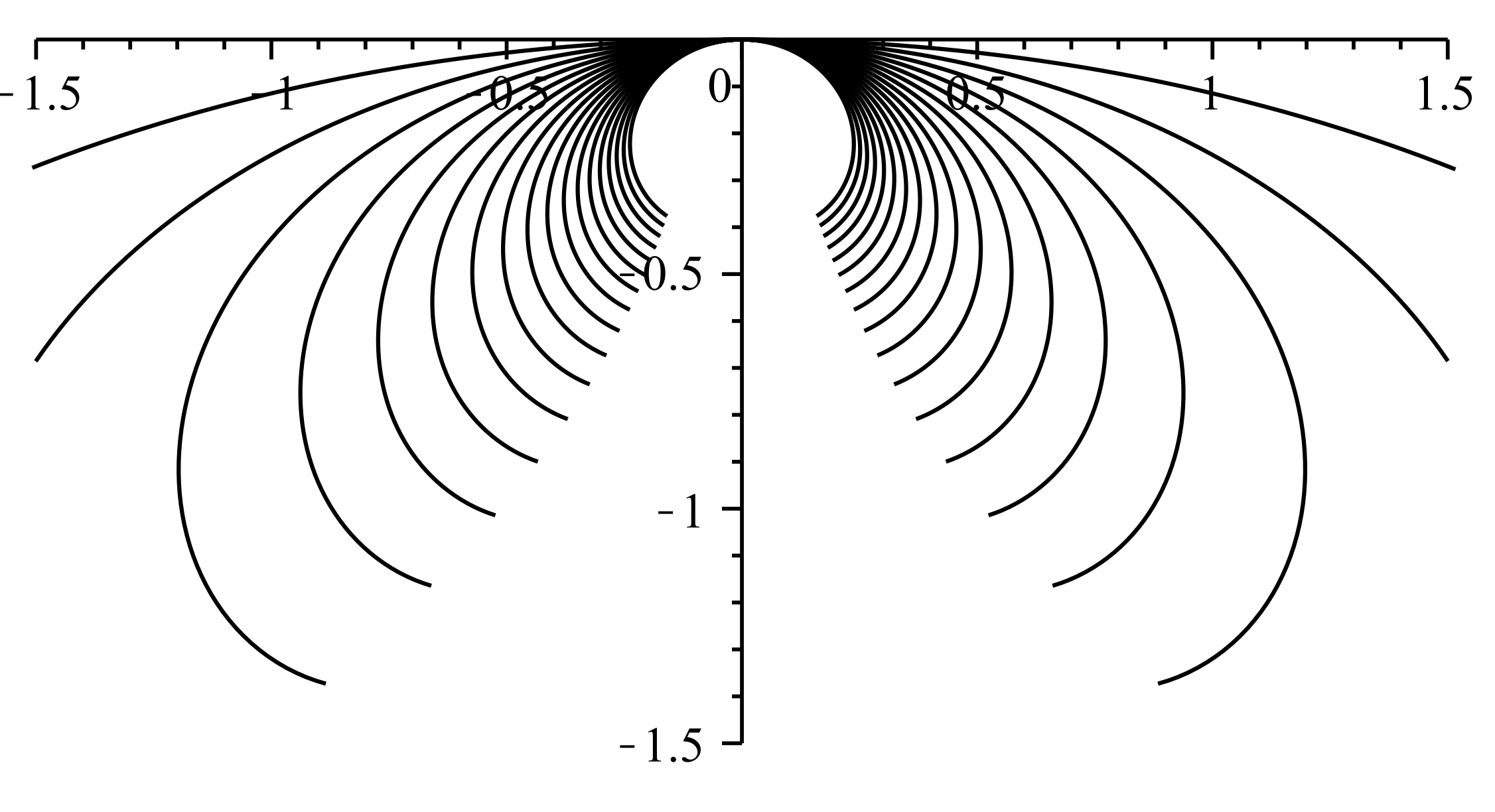}\\
  {Figure 6.6. Curves of $k_\varphi>1$}.
 \end{minipage}\\
  $\begin{cases}
          x(s)=2\,\arctan \left( {\dfrac {{{e}^{\sqrt {1-{c}^{2}}s}}-c}
          {\sqrt {1-{c}^{2}}}} \right) - cs,\\[0.4cm]
          y(s)=\ln\left( {{e}^{\sqrt{1-c^2}s}  + {e}^{-\sqrt{1-c^2}s}}- 2c\right).
    \end{cases}$&
  $\begin{cases}
             x(s)= 2\arctan\left(\sqrt{\dfrac{c-1}{c+1}}\tan\dfrac{\sqrt{c^2-1}}{2}s\right)-cs,\\[0.4cm]
             y(s)=-\ln\left(\dfrac{\tan^2\left(\dfrac{\sqrt{c^2-1}}{2}s\right)+1}
        {\dfrac{c-1}{c+1}\tan^2\left(\dfrac{\sqrt{c^2-1}}{2}s\right)+1}\right).
     \end{cases}$
\end{tabular}


\begin{thebibliography}{99}
\bibitem{Aa} M. A. S. Aarons, {\sl Mean curvature flow with a forcing term in Minkowski space}, Calc. Var. Partial Differential Equations 25 (2006), no. 2, 205-246.
\bibitem {caja} C. Carroll, A. Jacob, C. Quinn, R. Walters, {\sl The isoperimetric problem on planes with density}, Bull. Aust. Math. Soc. 78 (2008), no. 2, 177--197.
\bibitem{cami} A. Ca\~{n}ete, M. Miranda and D. Vittone, {\sl Some isoperimetric problems in planes with density}, J. Geo. Anal.,Vol. 20, No. 2, 243-290.
\bibitem {co2} I. Corwin, N. Hoffman, S. Hurder, V. Sesum, and Y. Xu, {\sl Differential geometry of manifolds with density}, Rose-Hulman Und. Math. J., 7 (1) (2006).
\bibitem {co3}   I. Corwin, F. Morgan, {\sl The Gauss-Bonnet formula on surfaces with densities}, Involve 4 (2011), no. 2, 199--202.
 \bibitem {dadu}   J. Dahlberg, A. Dubbs, E. Newkirk, H. Tran, {\sl Isoperimetric regions in the plane with density $r^p$}, New York J. Math. 16 (2010) 31-51.
  \bibitem {EcHui1}    K. Ecker, G. Huisken, {\sl Parabolic methods for the construction of spacelike slices of prescribed mean curvature in cosmological spacetimes}, Comm. Math. Phys. 135 (1991), no. 3, 595-613.
      \bibitem {Hiho}     D. T. Hieu, N. M. Hoang, {\sl Ruled minimal surfaces in $\Bbb R\sp 3$ with density $e\sp z$}, Pacific J. Math. 243 (2009), no. 2, 277-285.
    \bibitem {Hi}   D. T. Hieu, {\sl Some calibrated surfaces in manifolds with density}, J. Geom. Phys. 61 (2011), no. 8, 1625-1629.
 \bibitem {HuSi1} G. Huisken, C. Sinestrari, {\sl Mean curvature flow singularities for mean convex surfaces}, Calc. Var. PDE 8 (1999), no. 1, 1-14.
  \bibitem {HuSi2} G. Huisken, C. Sinestrari, {\sl  Convexity estimates for mean curvature flow and singularities of mean convex surfaces}, Acta Math. 183 (1999), no. 1, 45-70.
     \bibitem {Jian1}      H. Jian, H. Ju, Y. Liu, W. Sun, {\sl Symmetry of translating solutions to mean curvature flows}, Acta Math. Sci. Ser. B Engl. Ed. 30 (2010), no. 6, 2006-2016.
  \bibitem {Jian2}      H. Jian, H. Ju, Y. Liu, W. Sun,   {\sl Traveling fronts of curve flow with external force field}, Commun. Pure Appl. Anal. 9 (2010), no. 4, 975-986.
  \bibitem {Jian3}       H. Ju,  J. Lu, H. Jian,  {\sl Translating solutions to mean curvature flow with a forcing term in Minkowski space} Commun. Pure Appl. Anal. 9 (2010), no. 4, 963-973.
\bibitem{mamo} Q. Maurmann, F. Morgan, {\sl Isoperimetric comparison theorems for manifolds with density}, Calc. Var. PDE \textbf{36} (2009), No. 1, 1-5.
\bibitem {mo1} F. Morgan, {\sl Manifolds with density}, Notices Amer. Math. Soc., 52 (2005), 853-858.
\bibitem {mo3} F. Morgan, {\sl Myers' Theorem with density}, Kodai Math. J. 29 (2006), 454-460.
\bibitem {mo2} F. Morgan, ``Geometric Measure Theory: a Beginner's Guide'',
 Academic Press, fourth edition, 2009.
 \bibitem {mo5} F. Morgan, {\sl Manifolds with density and Perelman's proof of the Poincar\'{e} Conjecture},  Amer. Math. Monthly 116 (Feb., 2009), 134-142.
 \bibitem {NiTa}  H.   Ninomiya, M. Taniguchi, {\sl Traveling curved fronts of a mean curvature flow with constant driving force} Free boundary problems: theory and applications, I (Chiba, 1999), 206-221, GAKUTO Internat. Ser. Math. Sci. Appl., 13, Gakkotosho, Tokyo, 2000.
\bibitem {RCBM} C. Rosales, A. Ca\~{n}ete, V. Bayle and F. Morgan,
{\sl On the isoperimetric problem in Euclidean space with density}, Calc. Var. PDE \textbf{31} (2008), no. 1, 27-46.
\end{thebibliography}
\end{document}